\documentclass[11pt]{article}
\usepackage{latexsym}

\begin{document}

\newcommand{\ha}{H_{\alpha}}
\newcommand{\haa}{H_{\alpha}^A}
\newcommand{\hyp}{{\bf H^3}}
\newcommand{\xn}{\chi_n}
\newcommand{\sa}{S_{\alpha}}
\newcommand{\edel}{E_{\Delta}}

\newcommand{\calf}{{\cal F}}
\newcommand{\calg}{{\cal G}}
\newcommand{\pmf}{{\cal PMF}(S)}
\newcommand{\mf}{{\cal MF(S)}}
\newcommand{\calfn}{{\cal F}_n}
\newcommand{\calhn}{{\cal H}_n}
\newcommand{\calvn}{{\cal V}_n}
\newcommand{\tel}{{\cal T}_{el}(S)}
\newcommand{\boundary}{\partial_{\infty} \tel}
\newcommand{\calx}{{\cal X}}
\newcommand{\caly}{{\cal Y}}
\newcommand{\calz}{{\cal Z}}
\newcommand{\calh}{{\cal H}}
\newcommand{\pmfmin}{{\cal PMF}_{min}(S)}
\newcommand{\calv}{{\cal V}}
\newcommand{\ts}{{\cal T}(S)}
\newcommand{\cs}{{\cal C}(S)}
\newcommand{\csi}{{\cal C}_1(S)}
\newcommand{\thin}{Thin_{\alpha}}
\newcommand{\cso}{{\cal C}_0(S)}
\newcommand{\stuff}{p|_{[a_{i-1},a_i]}}
\newcommand{\callf}{{\it f}}
\newcommand{\fmin}{{\cal F}_{min}(S)}
\newcommand{\fs}{{\cal F}(S)}
\newcommand{\callx}{{\it x}}
\newcommand{\cally}{{\it y}}
\newcommand{\calan}{{\cal A}_n}
\newcommand{\calbn}{{\cal B}_n}
\newcommand{\hatx}{\hat{X}}

\newtheorem{thm}{Theorem}[section]
\newtheorem{prop}[thm]{Proposition}
\newtheorem{lemma}[thm]{Lemma}
\newtheorem{defin}[thm]{Definition}
\newtheorem{cor}[thm]{Corollary}

\begin{center}
{\Large The Boundary at Infinity of the Curve Complex and the Relative 
Teichm\"{u}ller Space}

\vspace{1mm}
Erica Klarreich

February 11, 1999
\end{center}

\begin{abstract}
In this paper we study the boundary at infinity of the curve complex $\cs$ of
a surface $S$ of finite type and the relative Teichm\"{u}ller space $\tel$
obtained from the Teichm\"{u}ller space by collapsing each region where a
simple closed curve is short to be a set of diameter 1.  $\cs$ and
$\tel$ are quasi-isometric, and Masur-Minsky have shown that $\cs$ and
$\tel$ are hyperbolic in the sense of Gromov.  We
show that the boundary at infinity of $\cs$ and $\tel$ is the space of
topological equivalence classes of minimal foliations on $S$.
\end{abstract}

\section{Introduction}

There is a strong but limited analogy between the geometry of
the Teichm\"{u}ller space ${\cal T}(S)$ of a surface $S$ and that
of hyperbolic spaces.  Teichm\"{u}ller space has many of the large-scale
qualities of hyperbolic space, and in fact the Teichm\"{u}ller space
of the torus is ${\bf H}^2$.  At one point it was generally believed
that the Teichm\"{u}ller metric was negatively curved; however,
Masur (\cite{masur}) showed that this is not so, apart from a few
exceptional cases.  Since then, Masur and Wolf (\cite{masur-wolf})
showed that ${\cal T}(S)$ is not even hyperbolic in the sense of
Gromov.

One way in which $\ts$ differs from hyperbolic space is that it does
not have a canonical compactification.  A Gromov hyperbolic space
has a boundary at infinity that is natural in the following two senses,
among others:  the boundary consists of all endpoints of quasigeodesic
rays up to equivalence (two rays are equivalent if they stay a bounded
distance from each other), and every isometry of the space
extends continuously to a homeomorphism of the boundary.
Teichm\"{u}ller space cannot be equipped with such a compactification but
rather gives rise to several compactifications, each with advantages 
and drawbacks.

Questions about the boundary of a hyperbolic space are interesting
for many reasons; one is that they tie in to questions of rigidity
of group actions by isometry on the space.  For example, in the
proof of Mostow's Rigidity Theorem, a key step in showing that
two hyperbolic structures on the same compact 3-manifold are
isometric is to show that a quasi-isometry between the two
structures lifts to a map of $\hyp$ that extends continuously
to $\partial_{\infty} \hyp$ (the Riemann sphere), and then to gain
some control over the map on $\partial_{\infty} \hyp$.  In another
instance, Sullivan's Rigidity Theorem gives geometric information
about a hyperbolic 3-manifold based on quasi-conformal 
information about its associated group action on $\partial_{\infty}
\hyp$.

Although Teichm\"{u}ller space is not hyperbolic, it is natural to be
interested in boundaries of Teichm\"{u}ller space, since they have a
strong connection to deformation spaces of hyperbolic 3-manifolds.
If $M$ is a compact 3-manifold, there is a well-known parametrization
of the space of geometrically finite hyperbolic structures on $int(M)$
by the Teichm\"{u}ller space of the boundary of $M$; see \cite{bers}.
One question is to understand the behavior of the hyperbolic structure
on $M$ as the Riemann surface structure on $\partial M$
``degenerates", that is, goes to infinity in the Teichm\"{u}ller space.
More generally, an important problem in the theory is to describe
all geometrically infinite hyperbolic structures on $M$; for this purpose
Thurston has introduced an invariant called the {\it ending lamination}
of $\partial M$, intended to play a similar role to that of the Teichm\"{u}ller
space of $\partial M$ in the geometrically finite setting.  Two important
boundaries of $\ts$ by Teichm\"{u}ller and Thurston involve compactifying
$\ts$ by the measured foliation space, or equivalently the measured
lamination space, which is related to but not the same as the space
of possible ending laminations on $S$.

Masur and Minsky (\cite{masur-minsky}) have shown that although
Teichm\"{u}ller space is not Gromov hyperbolic, it is {\it relatively
hyperbolic} with respect to a certain collection of closed subsets.
In this paper we describe the boundary at infinity of the relative Teichm\"{u}ller
space and a closely related object, the curve complex.  If $\alpha$ is
a homotopy class of simple closed curves on $S$, a surface of finite
type, let $Thin_{\alpha}$ denote the region of $\ts$ where the extremal
length of $\alpha$ is less than or equal to $\epsilon$, for some
fixed small $\epsilon > 0$.  These regions play a role somewhat
similar to that of horoballs in hyperbolic space; in fact for the torus,
these regions are actual horoballs in ${\bf H}^2$.  However, Minsky
has shown (\cite{min3}) that in general the geometry of each region
$Thin_{\alpha}$ is not hyperbolic, but rather has the large-scale
geometry of a product space with the sup metric.  On the other hand,
these regions are in a sense the only obstacle to hyperbolicity:
Masur and Minsky (\cite{masur-minsky}) have shown that Teichm\"{u}ller
space is relatively hyperbolic with respect to the family of regions
$\{ Thin_{\alpha} \}$.  In other words, the {\it electric Teichm\"{u}ller space}
$\tel$ obtained from $\ts$ by collapsing each region $Thin_{\alpha}$
to diameter 1 is Gromov hyperbolic (this collapsing is done by adding
a point for each set $Thin_{\alpha}$ that is distance $\frac{1}{2}$ from
each point in $Thin_{\alpha}$).

Since $\tel$ is Gromov hyperbolic, it can be equipped with a boundary
at infinity $\boundary$.  Our main result is the following:

\begin{thm}
\label{maina}
The boundary at infinity of $\tel$ is homeomorphic to the space of minimal
topological foliations on $S$.
\end{thm}

A foliation is minimal if no trajectory is a simple closed curve.
This space of minimal   foliations is exactly the space of possible
ending laminations (or foliations) on a surface $S$ that corresponds to
a geometrically infinite end of a hyperbolic manifold that has no parabolics 
(see \cite{ohshika}).  Here the topology on the space of minimal foliations is that
obtained from the measured foliation space by forgetting the measures.
This topology is Hausdorff (see the appendix), unlike the topology on
the full space of topological foliations; hence we may prove
Theorem \ref{maina} using sequential arguments to establish
continuity.

We will prove Theorem \ref{maina} by showing that the inclusion of $\ts$
in $\tel$ extends continuously to a portion of the Teichm\"{u}ller compactification
of $\ts$ by the projective measured foliation space $\pmf$:

\begin{thm}
\label{main}
The inclusion map from $\ts$ to $\tel$ extends continuously to the portion
$\pmfmin$ of $\pmf$ consisting of minimal foliations, to give a map
$\pi: \pmfmin \rightarrow \boundary$.  The map $\pi$ is surjective, and
$\pi(\calf)=\pi(\calg)$ if and only if $\calf$ and $\calg$ are topologically
equivalent.  Moreover, any sequence $\{ x_n \}$ in $\ts$ that converges
to a point in $\pmf \setminus \pmfmin$ cannot accumulate in the electric
space onto any portion of $\boundary$.
\end{thm}

If ${\cal F}_{min}(S)$ is the space of minimal topological foliations
on $S$, the map $\pi : \pmfmin \rightarrow \boundary$ descends to
a homeomorphism from ${\cal F}_{min}(S)$ to $\boundary$; hence
Theorem \ref{maina} is a consequence of Theorem \ref{main}.

\vspace{2mm}

Another space that we can associate to $S$ that has a close connection
to the electric Teichm\"{u}ller space is the {\it curve complex} ${\cal C}(S)$,
originally described by Harvey in \cite{harvey}.  
$\cs$ is a simplicial complex whose vertices are homotopy classes of
non-peripheral simple closed curves on $S$.  A collection of curves
forms a simplex if all the curves may be simultaneously realized so
that they are pairwise disjoint (when $S$ is the torus, once-punctured torus 
or four-punctured
sphere it is appropriate to make a slightly different definition;
see Section 4).  $\cs$ can be given a metric structure by assigning
to each simplex the geometry of a regular Euclidean simplex whose
edges have length 1.

In the construction of the electric Teichm\"{u}ller space, if the
value of $\epsilon$ used to define the sets $Thin_{\alpha}$ is
sufficiently small, then $Thin_{\alpha}$ and $Thin_{\beta}$
intersect exactly when $\alpha$ and $\beta$ have disjoint
realizations on $S$, that is, when the elements $\alpha$ and
$\beta$ of $\cs$ are connected by an edge.  Hence the 1-skeleton
$\csi$ of $\cs$ describes the intersection pattern of the sets $Thin_{\alpha}$;
$\csi$ is the nerve of the collection $\{ \thin \}$.  The
relationship between $\tel$ and $\cs$ is not purely topological.
Masur and Minsky have shown (\cite{masur-minsky}) that $\tel$
is quasi-isometric to $\csi$ and $\cs$.  This implies that $\csi$
and $\cs$ are also Gromov hyperbolic (although in the proof
if Masur and Minsky, the implication goes in the other direction).  
Two Gromov hyperbolic
spaces that are quasi-isometric have the same boundary at
infinity, so a consequence of Theorem \ref{maina} is the
following:

\begin{thm}
\label{boundary of curve complex}
The boundary at infinity of the curve complex $\cs$ is the space
of minimal   foliations on $S$.
\end{thm}

As with Teichm\"{u}ller space, the curve complex is important in the study
of hyperbolic 3-manifolds.  Let $M$ be a compact 3-manifold whose
interior admits a complete hyperbolic structure, and suppose $S$ is a component
of $\partial M$ that corresponds to a geometrically infinite end
$e$ of $M$.  Thurston, Bonahon and Canary
(\cite{thurston,bonahon,canary}) have shown that there is a sequence of
simple closed curves $\alpha_n \in \csi$ whose geodesic representatives
in $M$ ``exit the end $e$", that is, are contained in smaller and
smaller neighborhoods of $S$ in $M$.  Further, they showed that every
such sequence converges to a unique geodesic lamination (equivalently,
  foliation) on $S$.  In the case when the hyperbolic structure
on $int(M)$ has a uniform lower bound on injectivity radius, Minsky
has shown (\cite{min1,mina}) that the sequence $\{ \alpha_n \}$
is a quasigeodesic in $\csi$; a form of this was a key step in his proof of
the Ending Lamination Conjecture for such manifolds, giving quasi-isometric
control of the ends of $M$.

Since the sequence $\{ \alpha_n \}$ is a quasigeodesic in $\csi$, it must
converge to a point $\calf$ in the boundary at infinity of $\csi$, which
we have described as the space of minimal foliations (or laminations)
on $S$.  We will show that this description is natural, so that in
particular when the sequence $\{ \alpha_n \}$ in $\csi$ arises as
described above in the context of hyperbolic 3-manifolds, the boundary
point $\calf$ is the ending lamination.

\begin{thm}
\label{natural}
Let $\{ \alpha_n \}$ be a sequence of elements of $\csi$ that
converges to a foliation $\calf$ in the boundary at infinity of
$\cs$.  Then regarding the curves $\alpha_n$ as elements of the projective
measured foliation space $\pmf$, every accumulation point of $\{
\alpha_n \}$ in $\pmf$ is topologically equivalent to $\calf$.
\end{thm}

It is interesting to note that our description of the boundary
of $\cs$ ultimately does not depend on our original choice of a
Teichm\"{u}ller compactification for $\ts$, even though the
Teichm\"{u}ller boundary of $\ts$ depends heavily on an initial
choice of basepoint in $\ts$ (see Section 2 for more details).
Kerckhoff has shown (see \cite{kerckhoff}) that the action of
the modular group by isometry on $\ts$ does not extend
continuously to the Teichm\"{u}ller boundary; on the other hand,
the natural actions of the modular group on $\tel$ and $\cs$
do extend to the boundary at infinity, since this is true of
any action by isometry on a Gromov hyperbolic space.  Hence
the collapse used in the construction of $\tel$ essentially
``collapses" the discontinuity of the modular group action.

\vspace{2mm}

In Section 2 we will give an overview of some of the basic theory
of Teichm\"{u}ller space and quadratic differentials.  Section 3
contains the essential ideas of Gromov hyperbolicity that we will
need.  In Section 4 we discuss in more detail Masur and Minsky's
work on the electric Teichm\"{u}ller space and the curve complex,
and describe the quasi-isometry between them.  In Section 5
we establish some facts about convergence properties of
sequences of Teichm\"{u}ller geodesics, which are used in Section 6
to prove the main theorems.

\vspace{2mm}

\noindent {\bf Acknowledgments.}  The author would like to thank
Dick Canary and Yair Minsky for interesting conversations, and
Howard Masur for suggesting a portion of the argument for
Proposition \ref{segments}.


\section{Quadratic Differentials and the Teichm\"{u}ller Compactification
of Teichm\"{u}ller Space}

Let $S$ be a surface of finite genus and finitely many punctures.  The
Teichm\"{u}ller space ${\cal T}(S)$ is the space of all equivalence classes of
conformal structures of finite type on $S$, where two conformal structures
are equivalent if there is a conformal homeomorphism of one to the other
that is isotopic to the identity on $S$.  A conformal structure is of finite type
if every puncture has a neighborhood that is conformally equivalent to
a punctured disk.  The Teichm\"{u}ller distance between two points $\sigma$
and $\tau \in {\cal T}(S)$ is defined by
\[ d(\sigma,\tau) = \frac{1}{2} \log K(\sigma,\tau), \]
where $K(\sigma,\tau)$ is the minimal quasiconformal dilatation of any
homeomorphism from a representative of $\sigma$ to a representative of 
$\tau$ in the correct homotopy class.  The extremal map from $\sigma$ to $\tau$ may be constructed
explicity using quadratic differentials.

A holomorphic quadratic differential $q$ on a Riemann surface $\sigma$
is a tensor of the form $q(z) dz^2$ in local coordinates, where
$q(z)$ is holomorphic.  We define
\[ \| q \| = \int \int_S | q(z) | dx dy. \]
Let ${\cal DQ}(\sigma)$ denote the open unit ball in the space ${\cal Q}(\sigma)$ of
quadratic differentials on $\sigma$, and ${\cal SQ}(\sigma)$ the unit sphere.

Every $q \in {\cal DQ}(\sigma)$ determines a Beltrami differential $\| q \|
\frac{\overline{q}}{| q |}$ on $\sigma$, which in turn determines
a quasiconformal map from $\sigma$ to a new element $\tau_q$ of
${\cal T}(S)$; this map is the Teichm\"{u}ller extremal map between $\sigma$
and $\tau_q$.  The map that sends $q$ to $\tau_q$ is a homeomorphism,
giving an embedding of ${\cal T}(S)$ in ${\cal Q}(\sigma)$; ${\cal SQ}(\sigma)$ is the
boundary of ${\cal T}(S)$ in ${\cal Q}(\sigma)$, and ${\cal T}(S) \cup {\cal SQ}(\sigma)$ gives
a compactification of ${\cal T}(S)$ which we will denote $\overline{{\cal T}(S)}$, called the {\it Teichm\"{u}ller compactification} of ${\cal T}(S)$.

Any $q \in {\cal Q}(\sigma)$ determines a pair $\calh _q$ and $\calv _q$ of
measured foliations on $S$ called the horizontal and vertical foliations.
{\it Measured foliations} are equivalence classes of   foliations
of $S$ with 3- or higher-pronged saddle singularities, equipped with transverse measures; the equivalence is by
measure-preserving isotopy and Whitehead moves (that collapse
singularities).  We will denote the measured foliation space by
$\mf$ and the projectivized measured foliation space (obtained by
scaling the measures) by $\pmf$.  The horizontal and vertical foliations
associated to $q$ give a metric on $S$ in the conformal class of $\sigma$
that is Euclidean away from the singularities.  The map from ${\cal SQ}({\sigma})$
to $\pmf$ defined by sending $q$ to the projective class of its vertical
foliation is a homeomorphism, so that we may think of $\pmf$ as the 
boundary of ${\cal T}(S)$ (see \cite{hub-mas}).

A unit-norm quadratic differential $q$ on $\sigma$ determines a
directed geodesic line in ${\cal T}(S)$ as follows:  for $0 \leq k < 1$,
let $\sigma_k$ denote the element of ${\cal T}(S)$ determined by the
quasiconformal homeomorphism given by the quadratic differential
$k \cdot q$.  Geometrically, the extremal map from $\sigma$ to $\sigma_k$
is obtained by contracting the transverse measure of $\calh _q$ by a factor
of $K^{-{\frac{1}{2}}}$ and expanding the transverse measure of $\calv _q$ by $K^{\frac{1}{2}}$,
where $K= \frac{1+k}{1-k}$; note that the extremal map is $K$-quasiconformal.
The family $\{ \sigma_k : 0 \leq k < 1 \}$, when parametrized by
Teichm\"{u}ller arclength, gives a Teichm\"{u}ller geodesic ray; the family
$\{ \sigma_k : -1 < k < 1 \}$ determines a complete geodesic line.
Every ray and line through $\sigma$ is so determined.
We may think of the Teichm\"{u}ller ray $\{ \sigma_k : 0 \leq k < 1 \}$
as terminating at the boundary point $q \in {\cal SQ}(\sigma)$, or equivalently,
at the projective foliation $\calv _q \in \pmf$.  Similarly, every pair
of foliations in $\pmf$ that fills up $S$ (see the next section for
the definition of filling up) determines a geodesic line in ${\cal T}(S)$, for which
if $\tau \in L$ then the quadratic differential on $\tau$ that determines
$L$ has the two foliations as its horizontal and vertical foliations.

The compactification of $\overline{{\cal T}(S)}$ by endpoints of geodesic rays
depends in a fundamental way on the choice of basepoint $\sigma$ in
${\cal T}(S)$.  Kerckhoff has shown (see \cite{kerckhoff}) that there exist
projective foliations $\calf \in \pmf$ such that there are choices of
$\tau \in {\cal T}(S)$ for which the Teichm\"{u}ller ray from $\tau$ determined by
$\calf$ does not converge in $\overline{{\cal T}(S)}$ to $\calf$, but rather
accumulates onto a portion of $\pmf$ consisting of projective foliations
that are topologically equivalent but not measure equivalent to $\calf$.

\vspace{2mm}

\noindent {\bf Intersection number.}  If $\alpha$ is a simple closed
curve on $S$ then $\alpha$ determines a foliation on $S$ (which we
will also call $\alpha$) whose non-singular leaves are all freely
homotopic to $\alpha$.  The non-singular leaves form a cylinder,
and $S$ is obtained by gluing the boundary curves in some preassigned
manner.  There is a one-to-one correspondence between transverse measures
on $\alpha$ and positive real numbers:  each measure corresponds to the height
of the cylinder, that is, the minimal transverse measure of all arcs connecting
the two boundary curves of the cylinder.  If a measure has height $c$,
we will denote the measured foliation by $c \cdot \alpha$.  We define the
{\it intersection number} of the foliations $c \cdot \alpha$ and $k \cdot \beta$ by
\[ i(c \cdot \alpha, k \cdot \beta) = c k \cdot i(\alpha, \beta) \]
where the right-hand intersection number is just the geometric
intersection number of the simple closed curves $\alpha$ and $\beta$ 
(that is, the minimal number of crossings of any pair of representatives
of $\alpha$ and $\beta$.
Thurston has shown that the collection $\{ c \cdot \alpha : \alpha$ a simple
closed curve, $c \in {\bf R^+} \}$ is dense in $\mf$, and that the 
intersection number extends continuously to a function
$i: \mf \times \mf \rightarrow {\bf R}$ (see for instance \cite{aster}).

Note that although the intersection number of two {\it projective}
measured foliations is not well-defined, it still makes sense to
ask whether two projective measured foliations have zero or non-zero
intersection number.

A foliation $\calf$ is $minimal$ if no leaves of $\calf$ are simple
closed curves.  We say that two measured foliations are topologically
equivalent if the topological foliations obtained by forgetting
the measures are equivalent with respect to isotopy and Whitehead
moves that collapse the singularities.  Rees (\cite{rees}) has shown  the following:

\begin{prop}
\label{min int}
If $\calf$ is minimal then $i(\calf,\calg)=0$ if and only if
$\calf$ and $\calg$ are topologically equivalent.
\end{prop}

We say that two foliations $\calf$ and $\calg$ {\it fill up}
$S$ if for every foliation $\calh \in \mf$, $\calh$ has non-zero
intersection number with at least one of $\calf$ and $\calg$.
A consequence of Proposition \ref{min int} is that if $\calf$
is minimal, then whenever $\calg$ is not topologically equivalent
to $\calf$, $\calf$ and $\calg$ fill up $S$.

Let ${\cal F}_{min}(S)$ denote the space of minimal topological
foliations, with topology obtained from the space $\pmfmin$ of minimal
projective measured foliations by forgetting the measures.
Our goal is to show that ${\cal F}_{min}(S)$ is homeomorphic
to the boundary at infinity of the electric Teichm\"{u}ller
space.  We will use sequential arguments to show that certain
maps are continuous, so it is necessary to show the following,
whose proof can be found in the appendix:

\begin{prop}
\label{Hausdorff}
The space ${\cal F}_{min}(S)$ is Hausdorff and first countable.
\end{prop}

\noindent Note:  The entire space ${\cal F}(S)$ of topological
foliations on $S$ is not Hausdorff.  If $\alpha$ and $\beta$
are two distinct homotopy classes of simple closed curves that
can be realized disjointly on $S$, then regarded as topological
foliations, $\alpha$ and $\beta$ do not have disjoint
neighborhoods; every neighborhood of $\alpha$ or $\beta$ must
contain the topological foliation containing both $\alpha$
and $\beta$, and whose non-singular leaves are all homotopic
to $\alpha$ or $\beta$.

\vspace{2mm}

\noindent {\bf Extremal length.}  If $\gamma$ is a free homotopy
class of simple closed curves on $S$, an important conformal
invariant is the extremal length of $\gamma$, which is defined
as follows:

\begin{defin}
Let $\sigma \in {\cal T}(S)$, and let $\gamma$ be a homotopy class of simple
closed curves on $S$.  The extremal length of $\gamma$ on $\sigma$,
written $ext_{\sigma}(\gamma)$, is defined by
$\sup_{\rho} \frac{(l_{\rho}(\gamma))^2}{A_{\rho}}$,
where $\rho$ ranges over all metrics in the conformal class of
$\sigma$, $A_{\rho}$ denotes the area of $S$ with respect to $\rho$,
and $l_{\rho}(\gamma)$ is the infimum of the length of all
representatives of $\gamma$ with respect to $\rho$.
\end{defin}

\noindent Extremal length may be extended to scalar multiples
of simple closed curves by $ext(k \cdot \gamma) =
k^2 ext(\gamma)$, and extends continuously to the space of
measured foliations.

\vspace{1mm}

Our goal is to describe the boundary of the relative Teichm\"{u}ller
space $\tel$ obtained from ${\cal T}(S)$ by collapsing each of the
regions $Thin_{\gamma}$ of ${\cal T}(S)$ to be a set of bounded diameter,
where $Thin_{\gamma}$ is the region of ${\cal T}(S)$ where the simple
closed curve $\gamma$ has short extremal length.  We will need
the following lemma, which gives a connection between extremal
length and intersection number (see for instance \cite{min5}
Lemma 3.1-3.2 for a proof):

\begin{prop}
\label{ext}
Let $q$ be a quadratic differential with norm less than 1 on
$\sigma \in {\cal T}(S)$ with horizontal and vertical foliations $\calh$
and $\calv$, and let $\calf$ be a measured foliation on $S$.  
Then $ext_{\sigma}(\calf) \geq
(i(\calf,\calh))^2$; likewise $ext_{\sigma}(\calf) \geq
(i(\calf,\calv))^2$.
\end{prop}


\section{Gromov-hyperbolic spaces}

\vspace{3mm}

In this section we will present an overview of some of the basic theory
of Gromov-hyperbolic spaces.  References for the material in this section
are \cite{gdlh}, \cite{gr}, \cite{cdp} and \cite{can}.

Let $(\Delta,d)$ be a metric space.  If $\Delta$ is equipped with a basepoint
0, define the {\it Gromov product} $\langle x|y \rangle $ of the points $x$ and $y$ in
$\Delta$ to be
\[ \langle x|y \rangle =\langle x|y \rangle _0=\frac{1}{2} (d(x,0)+d(y,0)-d(x,y)).\]

\begin{defin}
  Let $\delta \geq 0$ be a real number.  The metric space $\Delta$ is
{\bf $\delta$-hyperbolic} if
\[ \langle x|y \rangle  \geq min( \langle x|y \rangle ,\langle y|z \rangle )- \delta \]
for every $x,y,z \in \Delta$ and for every choice of basepoint.
\end{defin}

We say that $\Delta$ is hyperbolic in the sense of Gromov if $\Delta$
is $\delta$-hyperbolic for some $\delta$.

\vspace{1mm}

A metric space $\Delta$ is {\it geodesic} if any two points in $\Delta$
can be joined by a geodesic segment (not necessarily unique).  If $x$
and $y$ are in $\Delta$ we write $[x,y]$, ambiguously, to denote some geodesic from
$x$ to $y$.

\vspace{2mm}

Heuristically, a $\delta$-hyperbolic space is ``tree-like'';  more precisely,
if we define an $\epsilon$-{\it narrow} geodesic polygon to be one such that every point on
each side of the polygon is at distance $\leq \epsilon$ from a point in
the union of the other sides, then we have

\begin{prop}
\label{narrow}
 In a geodesic $\delta$-hyperbolic metric space,
every n-sided polygon $(n \geq 3)$ is $4(n-2) \delta$-narrow.
\end{prop}

In a geodesic hyperbolic space, the Gromov product of two points $x$
and $y$ is roughly the distance from 0 to $[x,y]$; we have

\begin{prop}
\label{Gromov and dist from 0}
 Let $\Delta$ be a geodesic, $\delta$-hyperbolic
space and let $x,y \in \Delta$.  Then
\[ d(0,[x,y]) - 4\delta \leq \langle x|y \rangle  \leq d(0,[x,y]) \]
for every geodesic segment $[x,y]$.
\end{prop}

\vspace{3mm}

\noindent {\bf The boundary at infinity of a hyperbolic space.}  If $\Delta$ is a 
hyperbolic space, $\Delta$ can be equipped with a boundary
in a natural way.  We say that a sequence $\{ x_n \} $ of points in $\Delta$
{\it converges at infinity} if we have $\lim_{m,n \rightarrow \infty}
\langle x_m|x_n \rangle  = \infty$; note that this definition is independent of the
choice of basepoint, by Proposition \ref{Gromov and dist from 0}.  Given two sequences $\{ x_m \} $ and $\{ y_n \} $ that
converge at infinity, say that $\{ x_m \} $ and $\{ y_n \} $ are {\it equivalent}
if $\lim_{m,n \rightarrow \infty} \langle x_m|y_n \rangle  = \infty$.  Since $\Delta$
is hyperbolic, it is easily checked that this is an equivalence relation.  Define the {\it boundary at infinity}
$\partial_{\infty} \Delta$ of $\Delta$ to be the set of equivalence classes of
sequences that converge at infinity.  If $\xi \in \partial_{\infty} \Delta$ then
we say that a sequence of points in $\Delta$ {\it converges to $\xi$}
if the sequence belongs to the equivalence class $\xi$.  Write
$\overline{\Delta} = \Delta \cup \partial_{\infty} \Delta$.  When the space
$\Delta$ is a proper metric space, the boundary at infinity may
also be described as the set of equivalence classes of quasigeodesic
rays, where two rays are equivalent if they are a bounded
Hausdorff distance from each other.

\vspace{2mm}

\noindent {\bf Quasi-isometries and quasi-geodesics.}  Let $\Delta_0$ and $\Delta$
be two metric spaces.  Let $k \geq 1$ and $\mu \geq 0$ be real numbers. 
A quasi-isometry from $\Delta_0$ to $\Delta$ is a relation $R$ between
elements of $\Delta_0$ and $\Delta$ that has the coarse behavior of
an isometry.  Specifically, let $R$ relate every element of
$\Delta_0$ to some subset of $\Delta$ (so that we allow a given point in $\Delta_0$ to be
related to multiple points in $\Delta$).  We say that $R$ is a
$(k,\mu)$-{\it quasi-isometry} if for all $x_1$ and $x_2 \in \Delta_0$,
\[ \frac{1}{k} d(x_1,x_2) - \mu \leq d(y_1,y_2) \leq k d(x_1,x_2) + \mu \]
whenever $x_1 R y_1$ and $x_2 R y_2$.  Note that for a quasi-isometry,
given $x \in \Delta_0$ there is an upper bound to the diameter of the
set $\{ y \in \Delta | x R y \}$, that is independent of $x$.

We say that $R$ is a cobounded quasi-isometry if in addition, there is
some constant $L$ such that if $y \in \Delta$, $y$ is within $L$ of
some point that is related by $R$ to a point in $\Delta_0$.  If $R$
is a cobounded quasi-isometry then $R$ has a quasi-inverse, that is,
a relation $R^{\prime}$ that relates each element of $\Delta$ to
some subset of $\Delta_0$, with the following property:   there is
some constant $K$ for which if $x$ and $x^{\prime}$ are elements of 
$\Delta_0$ such that for some $y \in \Delta$, $x R y$ and 
$y R^{\prime} x^{\prime}$, then
$d(x,x^{\prime}) \leq K$.

\vspace{2mm}

A quasi-isometry between two $\delta$-hyperbolic spaces extends continuously
to the boundary, in the following sense:

\begin{thm}
\label{quasi extends}
 Let $\Delta_0$ and $\Delta$ be Gromov-hyperbolic, and let
$h:  \Delta_0 \rightarrow \Delta$ be a quasi-isometry.  For every
sequence $\{ x_n \} $ of points in $\Delta_0$ that converges to a point $\xi$ in
$\partial_{\infty} \Delta_0$, the sequence $\{ h(x_n) \} $ converges to a point in
$\partial_{\infty} \Delta$ that depends only on $\xi$, so that $h$ defines
a continuous map from $\partial_{\infty} \Delta_0$ to $\partial_{\infty} \Delta$.  The map
$h: \partial_{\infty} \Delta_0 \rightarrow \partial_{\infty} \Delta$ is injective.
\end{thm}

Theorem \ref{quasi extends} is, among other things, a key step 
in the proof of Mostow's
Rigidity Theorem.

\vspace{2mm}

If the metric on $\Delta$ is a path metric, 
a {\it $(k,\mu)$-quasigeodesic} is a
rectifiable path $p: I \rightarrow \Delta$, where $I$ is an
interval in {\bf R}, such that for all $s$ and $t$ in $I$,
\[ \frac{1}{k} l(p|_{[s,t]}) - \mu \leq d(p(s),p(t)) \leq k \cdot l(p|_{[s,t]}) + 
\mu . \]
Note that if a path $p: I \rightarrow \Delta$ is parametrized by arc length
then it is a quasigeodesic if and only if it is a quasi-isometry.

\vspace{1mm}

The behavior of quasigeodesics in the large is like that of actual
geodesics.  In particular, we have the following analogue of Proposition
\ref{Gromov and dist from 0}:

\begin{prop}
\label{segment distance}
Let $s:I \rightarrow \Delta$ be a quasigeodesic with endpoints $x$
and $y$.  Then there are constants $K$ and $C$ that only depend
on the quasigeodesic constants of $s$ and the hyperbolicity
constant of $\Delta$, such that
\[ \frac{1}{K} d(0,s(I)) - C \leq \langle x | y \rangle \leq
K d(0,s(I)) + C. \]
\end{prop}


\section{The Curve Complex and the Relative Hyperbolic Space}

\noindent {\bf The Curve Complex.}
If $S$ is an oriented surface of finite type, an important related
object is a simplicial complex called the {\it curve complex}.  Except 
in the cases
when $S$ is the torus, the once-punctured torus or a sphere with
4 or fewer punctures, we define the curve complex $\cs$
in the following way:  the vertices of $\cs$ are homotopy classes
of non-peripheral simple closed curves on $S$.  Two curves are
connected by an edge if they may be realized disjointly on $S$,
and in general a collection of curves spans a simplex if the curves
may be realized disjointly on $S$.

When $S$ is a sphere with 3 or fewer punctures, there are no
non-peripheral curves on $S$, so $\cs$ is empty.  When $S$
is the 4-punctured sphere, the torus, or the once-punctured
torus, there are non-peripheral simple closed curves on $S$,
but every pair of curves must intersect, so $\cs$ has no edges.  
For these three surfaces, a more interesting space to consider
is the complex in which two curves are connected by an edge
if they can be realized with the smallest intersection number
possible on $S$ (one for the tori; 2 for the sphere); we alter
the definition of $\cs$ in this way.  In these cases, $\cs$
is the Farey graph, which is well-understood (see for example
\cite{min2,min4}).

We give $\cs$ a metric structure by making every simplex
a regular Euclidean simplex whose edges have length 1.
The main result of \cite{masur-minsky} is the following:

\begin{thm}
\label{curve hyperbolic}
(Masur-Minsky)  $\cs$ is a $\delta$-hyperbolic space, where
$\delta$ depends on $S$.
\end{thm}

\noindent Note that $\cs$ is clearly quasi-isometric to its 1-skeleton
$\csi$, so that in particular $\csi$ is also Gromov-hyperbolic.

\vspace{2mm}

\noindent {\bf The Relative Teichm\"{u}ller Space.}
For a fixed $\epsilon > 0$, for each curve $\alpha \in \cso$
denote
\[ \thin = \{ \sigma \in \ts : ext_{\sigma}(\alpha) \leq \epsilon \} . \]
We will assume that $\epsilon$ has been chosen sufficiently
small that the collar lemma holds; in that case, a collection of
sets $Thin_{\alpha_1}, ... Thin_{\alpha_n}$ has non-empty
intersection if and only if $\alpha_1,...,\alpha_n$ can be realized
disjointly on $S$, that is, if $\alpha_1,...,\alpha_n$ form a simplex
in $\cs$.

We form the $relative$ or {\it electric Teichm\"{u}ller space} $\tel$
(following terminology of Farb \cite{farb}) by attaching a new point
$P_{\alpha}$ for each set $\thin$ and an interval of length $\frac{1}{2}$
from $P_{\alpha}$ to each point in $\thin$.  We give $\tel$ the
{\it electric metric} $d_{el}$ obtained from path length.

Masur and Minsky have shown the following:

\begin{thm}
\label{quasi-isometric}
\cite{masur-minsky}  $\tel$ is quasi-isometric to $\csi$.
\end{thm}

The quasi-isometry $R$ between $\csi$ and $\tel$ is defined as
follows:  if $\alpha$ is a curve in $\cso$, $\alpha$ is related to
the set $\thin$ (or equally well, to the ``added-on" point $P_{\alpha}$).
It is not difficult to see that the relation $R$ between $\cso$ and
$\tel$ is a quasi-isometry (see \cite{masur-minsky} for a proof).
$\cso$ is $\frac{1}{2}$-dense in $\csi$ (that is, every point in
$\csi$ is within $\frac{1}{2}$ of a point in $\cso$) and the collection
$\{ \thin \}$ is $D$-dense in $\tel$ for some $D$, so the relation
$R$ may easily be extended to be a cobounded quasi-isometry
from $\csi$ to $\tel$, making $\csi$ and $\tel$ quasi-isometric.

An immediate corollary of Theorem \ref{quasi-isometric} and
Theorem \ref{curve hyperbolic} is the following:

\begin{thm}
\label{electric hyperbolic}
\cite{masur-minsky}  The electric  Teichm\"{u}ller space $\tel$ is
hyperbolic in the sense of Gromov.
\end{thm}

\noindent We will use $\langle \cdot | \cdot \rangle_{el}$ to denote
the Gromov product on $\tel$.

\vspace{2mm}

\noindent {\bf Quasigeodesics in $\tel$.}  Since $\ts$ is contained
in $\tel$, each Teichm\"{u}ller geodesic is a path in $\tel$.  
Because certain portions of $\ts$ are collapsed to sets of bounded
diameter in $\tel$, a path whose Teichm\"{u}ller length is
very large may be contained in a subset of $\tel$ whose diameter
is small.  So to understand the geometry of these paths in $\tel$,
we introduce the notion of {\it arclength on the scale $c$},
after Masur-Minsky:
if $c>0$ and $p:[a,b] \rightarrow \tel$ is a path, we
define $l_c(p[a,b])=c \cdot n$ where $n$ is the smallest number
for which $[a,b]$ can be subdivided into $n$ closed subintervals
$J_1,...,J_n$ such that $diam_{\tel}(p(J_i)) \leq c$.

\vspace{1mm}

We will say that a path $p: [a,b] \rightarrow \tel$ in $\tel$ is an
electric quasigeodesic if for some $c>0$, $k \geq 1$ and $u > 0$ we have
\[ \frac{1}{k} l_c (p[s,t]) - \mu \leq
d_{el}(p(s),p(t)) \leq k \cdot l_c (p[s,t]) + \mu  \]
for all $s$ and $t$ in $[a,b]$ (note that the right-hand side
of the inequality is automatic).

Masur and Minsky have shown the following, which will be
important for understanding the boundary at infinity of
$\tel$:

\begin{thm}
\label{quasigeodesics}
\cite{masur-minsky}  Teichm\"{u}ller geodesics in $\ts$ are electric
quasigeodesics in $\tel$, with uniform quasigeodesic constants.
\end{thm}


\section{Convergence of sequences of Teichm\"{u}ller geodesics}

For the remainder of the paper we will assume that we have chosen
a basepoint $0 \in {\cal T}(S)$, giving an identification of ${\cal T}(S)$ with
the open unit ball of quadratic differentials on 0, and a
compactification $\overline{{\cal T}(S)}$ of ${\cal T}(S)$ by endpoints of
Teichm\"{u}ller geodesic rays from 0 (that is, by unit norm quadratic
differentials or equivalently, by projective measured foliations).

In view of Proposition \ref{segment distance} and the fact that
Teichm\"{u}ller geodesics are electric quasigeodesics, we can get some
control over the behavior of sequences going to infinity in the 
electric space if
we know the behavior of the Teichm\"{u}ller geodesic segments between
elements of the sequences.  The main fact we will need is the following:

\begin{prop}
\label{segments}

Let $\calf$ and $\calg$ be minimal foliations in $\pmf$.  Suppose $\{ x_n \}$
and $\{ y_n \}$ are sequences in ${\cal T}(S)$ that converge to $\calf$ and $\calg$,
respectively, and let $s_n$ denote the geodesic segment with endpoints
$x_n$ and $y_n$.  Then as $n \rightarrow \infty$, the sequence
$\{ s_n \}$ accumulates onto a set $s$ in $\overline{{\cal T}(S)}$ with
the following properties:

(1)  $s \cap {\cal T}(S)$ is a collection of geodesic lines whose horizontal
and vertical foliations are topologically equivalent to $\calf$ and $\calg$;
this collection is non-empty exactly when $\calf$ and $\calg$ fill up $S$
(that is, when $\calf$ and $\calg$ are not topologically equivalent).

(2)  $s \cap \partial {\cal T}(S)$ consists of foliations in $\pmf$ that
are topologically equivalent to $\calf$ or $\calg$.

\end{prop}

\noindent {\it Proof:}  We will begin by showing property {\it (2)}.
Let $\{ z_n \}$ be a sequence of points lying on the segments $s_n$
such that $z_n \rightarrow {\cal Z} \in \pmf$; we will show that
$\cal Z$ is topologically equivalent to either $\calf$ or $\calg$.

Suppose first that the $z_n$ lie over a compact region of moduli
space.  Then we claim that after dropping to a subsequence there is a sequence 
$\{ \alpha_n \}$ of distinct
simple closed curves on $S$ such that $ext_{z_n} (\alpha_n)$ is
bounded.  Since the $z_n$ lie over a compact region of moduli
space, there are elements $f_n$ of the mapping class group
that move the $z_n$ to some fixed compact region of Teichm\"{u}ller
space; since the $z_n$ are not contained in a compact region
of Teichm\"{u}ller space, we can drop to a subsequence so that
the maps $f_n$ are all distinct.  So, after dropping to a further
subsequence, there is some curve $\alpha$ on $S$ for which the
curves $\alpha_n=f_n^{-1}(\alpha)$ are all distinct; these curves
will have bounded extremal length on the surfaces $z_n$,
establishing the claim.  Now since $\pmf$ is compact, 
after dropping to a further
subsequence, the sequence $\{ \alpha_n \}$ converges in $\pmf$;
hence there exist constants $r_n$
such that the sequence $\{ r_n \alpha_n \}$
converges in $\mf$ to a foliation $\cal Z^{\prime}$, and since
the curves $\alpha_n$ are all distinct, we have $r_n \rightarrow
0$.  If instead
the $z_n$ do not lie over a compact region of moduli space then after
dropping to a subsequence
there is a sequence $\{ \alpha_n \}$ of (possibly non-distinct)
simple closed curves such that $ext_{z_n} (\alpha_n) \rightarrow
0$, and a sequence of bounded constants $r_n$ such that $r_n
\alpha_n$ converges to some $\cal Z^{\prime} \in \mf$.

Let $q_n$ denote the quadratic differential on the basepoint 0 that
is associated to $z_n$ by the identification of $\ts$ with ${\cal DQ}(0)$,
so that after dropping to a subsequence, $q_n \rightarrow q \in {\cal SQ}(0)$ whose vertical foliation
is $\cal Z$.  Let $\calfn$ denote the vertical foliation of $q_n$.
If we pull back $q_n$ by the Teichm\"{u}ller extremal map between
0 and $z_n$ to get a quadratic differential $\tilde{q}_n$ on $z_n$,
the vertical foliation of $\tilde{q}_n$ is $K_n^{1/2} \calfn$, where $K_n$
is the quasiconformality constant of the extremal map.  By Lemma
\ref{ext}, $ext_{z_n} (r_n \alpha_n) \geq (i(r_n \alpha_n, K_n^{1/2} \calfn ))^2$,
so $i(r_n \alpha_n, \calfn ) \rightarrow 0$ as $n \rightarrow \infty$.
So we have $i({\cal Z}^{\prime}, {\cal Z}) =0$.

On $z_n$, let $\phi_n$ denote the quadratic differential 
determining the segment $s_n$, and let ${\cal H}_n$ and
${\cal V}_n$ denote the horizontal and vertical foliations
associated to $\phi_n$ (so that as we move along $s_n$ in the
direction from $x_n$ to $y_n$, the transverse measure of
${\cal H}_n$ contracts and the transverse measure of ${\cal V}_n$
grows).
  Since $ext_{z_n}(\alpha_n) \geq
(i(\alpha_n,{\cal H}_n))^2$, we have $i(r_n \alpha_n, {\cal H}_n)
\rightarrow 0$; likewise $i(r_n \alpha_n, {\cal V}_n) 
\rightarrow 0$.  Let $a_n$ and $b_n$ be constants such that
after dropping to subsequences,
$a_n {\cal H}_n$ and $b_n {\cal V}_n$ converge to some 
$\cal H$ and $\cal V \in MF(S)$, respectively.  $\| \phi_n
\| = i(\calhn,\calvn)=1$, so since $i(\calh,\calv)$ must be
finite, the product $a_n b_n$ is bounded.  So we must have
at least one of the
sequences $\{ a_n\}$ and $\{ b_n\}$ bounded (say $\{ a_n \}$).
Then $i(r_n \alpha_n, a_n {\cal H}_n) \rightarrow 0$
as $n \rightarrow \infty$, so $i({\cal Z^{\prime}}, {\cal H})=0$.

Let $\tilde{\phi}_n$ denote the quadratic differential
on $x_n$ obtained by pulling back $\phi_n$ by the
Teichm\"{u}ller extremal map from $x_n$ to $z_n$.  Let $\tilde{\cal H}_n$
denote the horizontal foliation of $\tilde{\phi}_n$.  As we
move along $s_n$ from $z_n$ back to $x_n$ horizontal measure
grows, so we have $k_n \tilde{\cal H}_n = {\cal H}_n$ where
the constants $k_n$ are less than 1.  Following the argument
of the first paragraph of the proof, there is a sequence
$\{ \beta_n \}$ of simple closed curves on $S$ and a sequence
of bounded positive constants $t_n$ such that $ext_{x_n} (t_n \beta_n)
\rightarrow 0$ and $t_n \beta_n \rightarrow \calf^{\prime}$
where $i(\calf,\calf^{\prime}) = 0$ (so that $\calf^{\prime}$
is topologically equivalent to $\calf$, by minimality of
$\calf$).  This implies that $i(t_n \beta_n, \tilde{\cal H}_n)
\rightarrow 0$, so $i(t_n \beta_n, {\cal H}_n) \rightarrow 0$.
Taking limits, $i(\calf^{\prime}, {\cal H})=0$ so ${\cal H}$
is also topologically equivalent to $\calf$.  But we have already
shown that $i({\cal H}, {\cal Z}^{\prime}) = i({\cal Z}^{\prime},
{\cal Z}) = 0$, so by minimality $\cal Z$ is topologically
equivalent to $\calf$, establishing property {\it (2)}.

\vspace{1mm}

To show that $s \cap {\cal T}(S)$ consists of geodesic lines determined
by horizontal and vertical foliations topologically equivalent to
$\calf$ and $\calg$, suppose now that $\{ z_n \}$ is a sequence
of points in the segments $s_n$ such that $z_n \rightarrow z \in
{\cal T}(S)$.  Again, let $\phi_n$ denote the quadratic differential
on $z_n$ that determines the segment $s_n$, and let $\calhn$
and $\calvn$ denote the associated horizontal and vertical
foliations.  After descending to a subsequence, we can assume
that $q_n \rightarrow q$, a quadratic differential on $z$;
$\calhn$ and $\calvn$ will converge respectively to the horizontal
foliation $\cal H$ and vertical foliation $\cal V$ of $q$.
By arguments similar to those of the preceding paragraphs,
$\cal H$ and $\cal V$ are topologically equivalent to 
$\calf$ and $\calg$, respectively.  Now the segments $s_n$
all intersect a compact neighborhood of $z$, so since they form
an equicontinuous family of maps a subsequence must converge uniformly
on compact sets to the complete geodesic line containing $z$
determined by $q$.

When $\calf$ and $\calg$ are topologically equivalent,
it is impossible for any point in ${\cal T}(S)$ to support
a quadratic differential whose horizontal and 
vertical foliations are topologically equivalent to $\calf$
and $\calg$; hence when $\calf$ and $\calg$ are topologically
equivalent, $s \cap {\cal T}(S)$ must be empty.

\vspace{1mm}

It remains to show that when $\calf$ and $\calg$ fill $S$,
$s \cap {\cal T}(S)$ is nonempty.  We have shown that $s \cap \partial
{\cal T}(S)$ consists of foliations in $\pmf$ that are topologically
equivalent to $\calf$ or $\calg$.  The set of foliations in 
$\pmf$ topologically equivalent to $\calf$ is closed (likewise
for $\calg$), since if $\calf_n$ is a sequence of foliations
topologically equivalent to $\calf$ and $\calf_n \rightarrow
{\cal H} \in \pmf$ then by we have $i(\calf, {\cal H}) = 0$, 
so that $\cal H$ is
topologically equivalent to $\calf$.
So since $\calf$ and $\calg$ are not topologically equivalent,
$s \cap \partial {\cal T}(S)$ consists of at least two
connected components.  The segments $s_n$ are connected so
their accumulation set $s$ must be connected; hence
$s \cap {\cal T}(S)$ cannot be empty.  $\Box$ 

\vspace{2mm}

Note that in the course of the proof we have also shown the
following about sequences of segments whose endpoints converge
to foliations that are not minimal:

\begin{prop}
\label{nonminimal}
Let $x_n$ and $y_n$ be sequences in {\cal T}(S) converging to $\calf$ and
$\calg$ in $\pmf$, let $s_n$ be the geodesic segment with
endpoints $x_n$ and $y_n$, and let $s$ be the set of accumulation
points in $\overline{{\cal T}(S)}$ of the segments $s_n$.  Then the
only possible minimal foliations in $s \cap \pmf$ are those (if any)
that are topologically equivalent to $\calf$ or $\calg$.
\end{prop}

Using similar arguments we can prove the following about convergence
of Teichm\"{u}ller rays emanating from a common point (not
necessarily the chosen basepoint 0 in ${\cal T}(S)$):

\begin{prop}
\label{rays}
Let $z$ be a fixed point in ${\cal T}(S)$, let $z_n$ be a sequence of
points in ${\cal T}(S)$ that converge to ${\calz} \in \pmf$, and let
$r_n$ be the geodesic segment from $z$ to $z_n$.  After
descending to a subsequence, the segments $r_n$ converge 
uniformly on compact sets to a geodesic ray $r$ with vertical
foliation $\calv$, such that $i(\calz,\calv)=0$.
\end{prop}

\noindent {\it Proof:}  Assume that the segments $r_n$ are
paths parametrized by arclength, and extend the $r_n$ to maps
$r_n : {\bf R} \rightarrow \ts$ by setting $r_n(t)=z_n$ for
all $t \geq d(z,z_n)$.  The family $\{ r_n \}$ is equicontinuous,
so by Ascoli's theorem, after dropping to a subsequence the
maps $r_n$ converge uniformly on compact sets to a map
$r: {\bf R} \rightarrow \ts$, which is necessarily a geodesic
ray emanating from $z$.

Let $\calv$ be the vertical foliation of $r$.  We wish to show that
$i(\calz,\calv)=0$.  Let $\phi_n$ be the quadratic differential
on $z$ determining the segment $r_n$, and let $\calv_n$ be the
vertical foliation of $\phi_n$, so that $\calv_n \rightarrow
\calv$.  Then $ext_z \calv_n \rightarrow ext_z \calv$, so 
$ext_z \calv_n$ is bounded.  Now $d(z,z_n)=
\frac{1}{2}\log(\frac{ext_z \calv_n}{ext_{z_n} \calv_n})$ (see \cite{kerckhoff}), so
since $d(z,z_n) \rightarrow \infty$, $ext_{z_n} \calv_n
\rightarrow 0$ as $n \rightarrow \infty$.  Now the argument
of the third paragraph of the proof of Proposition \ref{segments}
(changing the $r_n \alpha_n$ to $\calv_n$)  shows that
$i(\calz,\calv)=0$.  $\Box$


\section{The boundary of the relative Teichm\"{u}ller space}

As a start to proving Theorem \ref{main} we will prove the following,
which shows that minimal foliations in $\pmf$ are an infinite
electric distance from any point in ${\cal T}(S)$.

\begin{prop}
\label{minimal is infinite distance}
Let $\calf \in \pmf$ be minimal and let $\{ z_n \}$ be
a sequence of points in ${\cal T}(S)$ that converges to $\calf$.  Then
$d_{el}(0,z_n) \rightarrow \infty$ as $n \rightarrow \infty$.
\end{prop}

\noindent {\it Proof:}  Suppose that $d_{el}(0,z_n)$ does
not go to infinity.  Then after dropping to a subsequence
we may assume that the $z_n$ lie in a bounded electric neighborhood
of 0.  As in the proof of Proposition \ref{segments}, we can
construct a sequence of curves $\alpha_n$ such that the values
$ext_{z_n}(\alpha_n)$ are bounded, and such that for some
bounded constants $r_n$, the sequence $r_n \alpha_n$ converges
in $\mf$ to a foliation $\calf_0$ such that $i(\calf,\calf_0)=0$.
Now since $\alpha_n$ has bounded extremal length on $z_n$,
we have that $z_n$ lies in a bounded neighborhood of $Thin_{\alpha_n}$,
so the values $d_{el}(0,Thin_{\alpha_n})$ are bounded.  So the curves
$\alpha_n$, regarded as elements of the curve complex, are a
bounded distance (say $M$) from some fixed curve $\alpha$.
Now for each $\alpha_n$ we can construct a chain of curves
$\alpha_{n,0},...,\alpha_{n,M}$ such that $\alpha_{n,0}=\alpha_n$
and $\alpha_{n,M}=\alpha$, and for all $i$, $d(\alpha_{n,i},\alpha_{n+1,i})
=1$.  So $\alpha_{n,i}$ and $\alpha_{n,i+1}$ are disjoint, or
in other words, $i(\alpha_{n,i},\alpha_{n,i+1})=0$.
After dropping to subsequences, for each fixed $i$, the 
sequence $\alpha_{n,i}$ converges (after bounded rescaling)
to a measured foliation $\calf_i$, and for all $i$ we have
$i(\calf_i,\calf_{i+1})=0$.  Since $\calf$ is minimal this
implies that all the foliations $\calf_i$ are topologically
equivalent to $\calf$.  But $\calf_M=\alpha$, which gives
a contradition.  $\Box$

\vspace{2mm}

The proof of Theorem \ref{main} will be divided into the next
three propositions.  We begin by showing that we have a well-
defined, continuous map from $\pmfmin$ to $\boundary$.

\begin{prop}
\label{extend}
The inclusion map from ${\cal T}(S)$ to $\tel$ extends continuously to
the portion $\pmfmin$ of $\pmf$ consisting of minimal foliations,
to give a map $\pi: \pmfmin \rightarrow \boundary$.
\end{prop}

\noindent {\it Proof:}  Let $\calf \in \pmfmin$.  We must
show that every sequence $\{ z_n \}$ in ${\cal T}(S)$ converging
to $\calf$, considered as
a sequence in $\tel$, converges to a
unique point in $\boundary$.  So suppose that there is a sequence
$\{ z_n \} \rightarrow \calf$ that does not converge to a point in 
$\boundary$.  Then there are subsequences $\{ x_n \}$ and
$\{ y_n \}$ of $\{ z_n \}$ such that $\langle x_n |
y_n \rangle_{el}$ is bounded.  Let $s_n$ denote the
Teichm\"{u}ller geodesic segment between $x_n$ and $y_n$.  Since the
segments $s_n$ are electric quasigeodesics with uniform
quasigeodesic constants, by Proposition \ref{segment distance}
there is a point $p_n$ on each $s_n$ that is a bounded electric
distance from 0.  By Proposition \ref{segments}, the points
$p_n$ converge to a foliation in $\pmfmin$ that is
topologically equivalent to $\calf$.  But then according to
Proposition \ref{minimal is infinite distance},
$d_{el}(0,p_n)$ must go to infinity as $n \rightarrow \infty$.
This gives a contradiction.  $\Box$

\vspace{2mm}

We now show that the non-injectivity of the map $\pi :
\pmfmin \rightarrow \boundary$ is limited to identifying 
foliations that are topologically equivalent but not measure
equivalent.

\begin{prop}
\label{identifications}
Let $\calf$ and $\calg$ be minimal foliations in $\pmf$.  Then
$\pi (\calf) = \pi (\calg)$ if and only if $\calf$ and $\calg$
are topologically equivalent.
\end{prop}

\noindent {\it Proof:}  Suppose first that $\calf$ and $\calg$
are topologically equivalent, and suppose that $\pi(\calf) \neq
\pi(\calg)$.  Then the same argument as in the proof of
Proposition \ref{extend} would give a sequence of points
$\{ p_n \}$ that are a bounded electric distance from 0
and that converge to a minimal foliation in $\pmf$; but this
is impossible by Proposition \ref{minimal is infinite distance}.
Hence when $\calf$ and $\calg$ are topologically equivalent,
$\pi(\calf)=\pi(\calg)$.

Now suppose that $\calf$ and $\calg$ are not topologically 
equivalent, and let $\{ x_n \}$ and $\{ y_n \}$ be sequences
in ${\cal T}(S)$ converging to $\calf$ and $\calg$, respectively.
We will show that $\{x_n \}$ and $\{ y_n \}$ do not
converge to the same point in $\boundary$, by showing that
we can drop to subsequences so that $\langle x_n |
y_n \rangle_{el}$ is bounded as $n \rightarrow \infty$.
Let $s_n$ denote the Teichm\"{u}ller geodesic segment with endpoints
$x_n$ and $y_n$.  Since $\calf$ and $\calg$ are not topologically
equivalent, by Proposition \ref{segments} we can drop to a
subsequence so that the $s_n$ converge uniformly on compact
sets to a Teichm\"{u}ller geodesic line $L$.  Choose a point
$p \in L$, and a sequence $p_n \in s_n$ converging to $p$.
Then as $n \rightarrow \infty$, $d(0,p_n)$ is bounded, 
hence $d_{el}(0,\pi(p_n))$ is also bounded.  So by Proposition
\ref{segment distance}, $\langle x_n | y_n \rangle_{el}$
is bounded as $n \rightarrow \infty$.  Thus $\pi(\calf) \neq
\pi(\calg)$.  $\Box$

\vspace{2mm}

The following proposition completes the proof of Theorem \ref{main}.

\begin{prop}
\label{surjective}
The map $\pi: \pmfmin \rightarrow \boundary$
is surjective.  Moreover, if $\{ x_n \}$ is a sequence in $\ts$ that
converges to a non-minimal foliation in $\pmf$ then no subsequence
of $\{ x_n \}$ converges in the electric space $\tel$ to a point in
$\boundary$.
\end{prop}

\noindent {\it Proof:}  Let $\calx \in \boundary$, and let
$x_n$ be a sequence in $\tel$ that converges to $\calx$; without
loss of generality we may assume that each $x_n$ lies in
${\cal T}(S)$, since if $x_n$ is one of the added-on points in the construction
of $\tel$ then we may replace $x_n$ by a point in ${\cal T}(S)$ that
is distance $\frac{1}{2}$ from $x_n$, without changing the convergence
properties of the sequence $\{ x_n \}$.  We will show that a
subsequence of $\{ x_n \}$ converges to a minimal foliation
$\calf \in \pmf$; then $\pi (\calf) = \calx$.

Since $\overline{{\cal T}(S)}$ is compact, after dropping to a subsequence,
$\{ x_n \}$ converges to some $\calf \in \pmf$.  Suppose $\calf$ 
is not minimal.  We will show that for some $B < \infty$, for
each $x_n$ there are infinitely many $x_m$ such that 
$\langle x_n|x_m \rangle_{el} < B$; this would contradict
convergence in $\overline{\tel}$ of the sequence $\{x_n \}$.  
Fix $x_n$, and let
$r_{mn}$ denote the geodesic segment with endpoints $x_n$ and $x_m$.
By Proposition \ref{rays}, a subsequence of the $r_{mn}$ (which
we will again call $r_{mn}$) converges uniformly on compact sets
to a geodesic ray $r_n$.  Let $\calhn$ denote the horizontal 
foliation of $r_n$; by Proposition \ref{rays} we have $i(\calf,
\calhn)=0$.  The foliations $\calf$ and $\calhn$ are not minimal,
so each one
contains a simple closed curve, which we will denote $\alpha$
and $\gamma_n$, respectively.  Now we have $i(\alpha,
\gamma_n)=0$, so that in the curve complex,
the distance from $\alpha$ to $\gamma_n$ is at most 1; hence
the electric distance from $Thin_{\alpha}$ to $Thin_{\gamma_n}$
is bounded independent of $n$, since the curve complex is
quasi-isometric to $\tel$.

The simple closed curve $\gamma_n$ contained in $\calhn$ may
be chosen so that as $t \rightarrow \infty$, $ext_{r_n(t)}
\gamma_n \rightarrow 0$ (see \cite{min}, Lemma 8.3).
So for all sufficiently large $t$, $r_n(t)$
belongs to $Thin_{\gamma_n}$.  Since the rays $r_{mn}$ converge
to $r_n$ uniformly on compact sets, for all $m$ sufficiently
large there is a point $p_{mn}$ on $r_{mn}$ that lies in 
$Thin_{\gamma_n}$.  Now we have
\[ d_{el}(0,p_{mn}) \leq d_{el}(0,Thin_{\gamma_n})+1
\leq d_{el}(0,Thin_{\alpha})+d_{el}(Thin_{\alpha},Thin{\gamma_n})+2 \]
since each thin set has diameter 1 (here the electric distance
between two sets $S_1$ and $S_2$ means the smallest distance
between any pair of points in $S_1$ and $S_2$, respectively).
Note that the right-hand side of the inequality does not depend
on $n$ or $m$ since $d_{el}(Thin_{\alpha},Thin_{\gamma_n})$
is bounded independent of $n$.  Now by Proposition \ref{segment distance}
we have that for all $m$ sufficiently large,
$\langle x_n | x_m \rangle_{el}$ is bounded, and
the bound does not depend on $n$ or $m$; this contradicts
the fact that the sequence $\{ x_n \}$ converges to a
point in the boundary at infinity of $\tel$, so our assumption
that $\calf$ is not minimal must be false.  Hence $\calf$ is
minimal, and we have $\pi(\calf) = \calx$.  $\Box$

\vspace{2mm}

Note that given a nonminimal foliation $\calf$, there are sequences in
${\cal T}(S)$ converging to $\calf$ whose electric distance from 0
goes to infinity; however, no subsequences of these will converge to a point
in $\boundary$, so that in particular $\tel \cup \boundary$
is not compact.  It is simple to construct such sequences:  the
minimal foliations are dense in $\pmf$ (see for instance
\cite{aster}), so there is a sequence $\{ \calf_n \}$ of
minimal foliations that converges to $\calf$.  By Proposition 
\ref{minimal is infinite distance}, for every $M>0$,
each $\calf_n$ has a neighborhood whose points are all at least $M$ 
from 0 in the electric metric; hence we
may easily choose a sequence $\{ p_n \}$ of points contained in small
neighborhoods of the foliations $\calf_n$, such that $\{ p_n
\}$ converges to $\calf$ and $d_{el}(0,p_n) \rightarrow \infty$.

\vspace{2mm}

If $\calf$ is a foliation in $\pmf$, let $\tau(\calf)$ denote
the equivalence class of foliations in $\pmf$ that are topologically
equivalent to $\calf$.  We have shown that the boundary at infinity
of  $\tel$ and $C(S)$ can be identified with topological equivalence classes of minimal
foliations.  In spite of the fact that the Teichm\"{u}ller 
compactification of ${\cal T}(S)$ by $\pmf$ depends heavily on the choice 
of basepoint, the arguments we have given show that the description we have 
obtained of the boundary of  $C(S)$ is natural:

\vspace{2mm}

\noindent {\bf Theorem 1.4}  {\it Let $\{ \alpha_n \}$ be a sequence of elements of $\csi$ that
converges to a foliation $\calf$ in the boundary at infinity of
$\cs$.  Then regarding the curves $\alpha_n$ as elements of the projective
measured foliation space $\pmf$, every accumulation point of $\{
\alpha_n \}$ in $\pmf$ is topologically equivalent to $\calf$.}


\section{Appendix}

In order to use sequential arguments to prove the continuity results
of the main theorems, it is necessary to understand the point-set
topology of $\fmin$, the space of minimal topological foliations
on $S$.  This is particularly important in light of the fact that
the entire space $\fs$ of topological foliations, with the topology
induced from $\pmf$ by forgetting the measures, is not Hausdorff.
We will begin with the following:

\begin{prop}
The measure-forgetting quotient map $p: \pmfmin \rightarrow \fmin$
is a closed map, and the pre-image of any point of $\fmin$ is compact.
\end{prop}

\noindent {\it Proof:}  To show that $p$ is a closed map, let 
$K \subset \pmfmin$ be a closed set.  Then
we claim that the set $p^{-1}(p(K))$ is closed; this will imply
that $p(K)$ is closed.  So, let $\{ x_n \}$ be s sequence in
$p^{-1}(p(K))$ that converges to a point $x$ in $\pmfmin$;
we must show that $x \in p^{-1}(p(K))$.
There is a sequence of $y_n \in K$ such that $p(x_n)=p(y_n)$.
Since $\pmf$ is compact, after dropping to a subsequence
we may assume that $y_n \rightarrow y \in \pmf$.  Now since
$p(x_n)=p(y_n)$, we have that $x_n$ and $y_n$ are topologically
equivalent, which implies that $i(x_n,y_n)=0$.  Hence
$i(x,y)=0$, so since $x$ is minimal, $x$ and $y$ are topologically
equivalent by Proposition \ref{min int}, so that $p(x)=p(y)$.
We now know $y$ to be in $\pmfmin$, so since $K$ is closed in
$\pmfmin$, we have $y \in K$.  This in turn implies that
$x \in p^{-1}(p(K))$, so $p^{-1}(p(K))$ is closed.

To show that the pre-image of any point is compact, let $z$ be
a point in $\fmin$ and let $Z=p^{-1}(z)$.  Let $\{ x_n \}$
be a sequence of points in $Z$; since $\pmf$ is compact,
after dropping to a subsequence we may assume that $x_n$
converges to some $x \in \pmf$.  Let $y$ be a fixed point
in $Z$.  Then the set $Z$ is the set of all foliations
in $\pmf$ that are topologically equivalent to $y$.  Hence
$i(y,x_n)=0$ for all $n$, so $i(y,x)=0$.  Thus by minimality,
$x$ is topologically equivalent to $y$, so $x \in Z$.
So $Z$ is compact.  $\Box$

\vspace{2mm}

The space $\pmf$ is metrizable and normal, since it is a
topological sphere; hence so is $\pmfmin \subset \pmf$.
The following proposition will establish in particular
that $\fmin$ is first countable and Hausdorff, which are
exactly the properties needed in order for sequential
arguments to prove continuity:

\begin{prop}
Let $X$ be a metric space that is normal, and let
$p: X \rightarrow \hatx$ be a quotient map that is a closed
map, and such that the pre-image of any point of $\hatx$ is compact.
Then the quotient topology on $\hatx$ is first countable
and normal.
\end{prop}

\noindent {\it Proof:}  We will show first that $\hatx$ is
normal.  Let $S$ and $T$ be disjoint closed sets in $\hatx$; we must
show that $S$ and $T$ have disjoint neighborhoods.  The sets
$p^{-1}(S)$ and $p^{-1}(T)$ are closed and disjoint in $X$,
so since $X$ is normal there are disjoint open sets $U$ and
$V$ such that $p^{-1}(S) \subset U$ and $p^{-1}(T) \subset V$.
Then $X - U$ and $X - V$ are closed, so
$p(X - U)$ and $p(X - V)$ are closed since
$p$ is a closed map.  Now $S$ has empty intersection with
$p(X - U)$, so $\hatx - p(X - U)$
is a neighborhood of $S$; likewise $\hatx - p(X - V)$
is a neighborhood of $T$.  It is easily checked that the sets
$\hatx - p(X - U)$ and $\hatx -
p(X - V)$ are disjoint, which establishes normality.

To show that $\hatx$ is first countable, let $z \in \hatx$;
we must define a countable neighborhood basis around $z$.
Let $Z = p^{-1}(z)$, and let $U_n$ be the open neighborhood
around $Z$ of radius $\frac{1}{n}$.  Let $V_n = p(U_n)$.
If $V$ is any neighborhood of $z$, then $p^{-1}(V)$ is
a neighborhood of the set $Z$, so since by assumption $Z$
is compact, $p^{-1}(V)$ must contain one of the sets $U_n$;
hence $V$ must contain one of the sets $V_n$.  So we will
be done if we can show that every $V_n$ contains a 
neighborhood of $z$.  In $X$, let $W_n = p^{-1}(V_n)$;
note that $U_n \subset int(W_n)$.  Let $S_n =
X - int(W_n)$, so that $S_n \cap U_n = \emptyset$.
The set $p(S_n)$ is closed in $\hatx$ since $p$ is
a closed map, so by normality of $\hatx$, there is some
neighborhood $T_n$ of $x$ disjoint from $p(S_n)$.
Now $p^{-1}(T_n) \subset W_n$, so $T_n \subset p(W_n)=V_n$.
Hence the sets $T_n$ form a local basis of neighborhoods
of $z$.  $\Box$




\vspace{3mm}

\noindent Department of Mathematics, University of Michigan, East Hall, Ann Arbor, MI 48109-1109; klarreic@math.lsa.umich.edu

\end{document}